\documentclass[a4paper,11pt]{amsart}
\usepackage{amsmath}
\usepackage{amssymb}
\usepackage{amsfonts}
\usepackage{amsthm}
\usepackage[hidelinks]{hyperref}
\usepackage[margin=1.2in]{geometry}

\usepackage[
backend=bibtex,
style=numeric,
sorting=nyt
]{biblatex}
\title{The weak Gram law for Hecke \( L \)-functions}
\author{Sebastian Weishäupl}
\address{Department of Mathematics, Würzburg University, 97074 Würzburg, Germany}
\email{sebastian.weishaeupl@mathematik.uni-wuerzburg.de}
\date{\today}
\subjclass[2010]{11M41, 11F66}
\keywords{Gram's law, Gram points, Hecke \( L \)-functions}

\newcommand{\eps}{\varepsilon}
\newcommand{\R}{\mathbb{R}}
\newcommand{\C}{\mathbb{C}}

\newcommand{\Z}{\mathbb{Z}}

\DeclareMathOperator{\Real}{Re}
\DeclareMathOperator{\Imag}{Im}
\DeclareMathOperator{\SL}{SL}

\DeclareMathOperator{\res}{res}

\newtheorem{theorem}{Theorem}
\newtheorem{lemma}[theorem]{Lemma}
\newtheorem{corollary}[theorem]{Corollary}

\theoremstyle{definition}

\begin{document}
\begin{abstract}
	We generalize a theorem by Titchmarsh about the mean value of Hardy's \( Z \)-function at the Gram points to the Hecke \( L \)-functions, which in turn implies the weak Gram law for them. Instead of proceeding analogously to Titchmarsh with an approximate functional equation we employ a different method using contour integration.
\end{abstract}
\maketitle
\section{Introduction}
	In the theory of the Riemann zeta-function the weak Gram law makes a statement about the distribution of the zeros of \( \zeta(s) \) on the critical line. To describe this statement we need several definitions. Starting from the	functional equation of \( \zeta(s) \) in its asymmetric form
	\begin{equation}\label{zetaFE}
		\zeta(s) = \Delta_\zeta(s) \zeta(1-s)\quad\;\;\text{with} \;\; \Delta_\zeta(s) = \pi^{s - \frac{1}{2}} \frac{\Gamma(\frac{1-s}{2})}{\Gamma(\frac{s}{2})}
	\end{equation}
	we define the function \( \vartheta_\zeta(t) \) as the continuous branch of the argument of \( \Delta_\zeta(\frac{1}{2}+it)^{-1/2} \) for \( t \in \R \) with \( \vartheta_\zeta(0) = 0 \). Hardy's \( Z \)-function is then defined by
	\[
		Z_\zeta(t) := e^{i \vartheta_\zeta(t)} \zeta\Big(\frac{1}{2}+it\Big).
	\]
	From the functional equation (\ref{zetaFE}) it follows that \( Z_\zeta(t) \) is real-valued. Furthermore the ordinates of the zeros of \( \zeta(s) \) on the critical line coincide with the zeros of \( Z_\zeta(t) \).
	
	The function \( \vartheta_\zeta(t) \) increases monotonically for \( t \geq 7 \) and grows arbitrarily large. This allows us to define the Gram points \( t_v \) as the unique solutions of \( \vartheta_\zeta(t_v) = v \pi \) for integers \( v \geq -1 \). These points were first studied by J\o rgen Pedersen Gram \cite{gram} in 1903 in the context of numerical computations of the zeros of \( \zeta(s) \). He observed that the Gram points and the ordinates of the zeros of \( \zeta(s) \) on the critical line (i.e. the zeros of \( Z_\zeta(t) \)) seem to alternate. Hutchinson \cite{hutchinson} called this phenomenon Gram's law and showed that it first fails in the interval \( [t_{125}, t_{126}] \), because it contains no zero of \( Z_\zeta(t) \). Titchmarsh \cite{titchmarsh} (see also \cite[\S10.6]{titchmarshZeta}) proved a mean value result for \( Z_\zeta(t) \) at the Gram points,  namely that we have for any fixed integer \( M \geq -1 / 2 \), as \( N \to \infty \),
		\begin{align*}
			\sum_{v=M}^N Z_\zeta(t_{2v}) &= 2N + O\big(N^{\frac{3}{4}} \log(N)^\frac{3}{4}\big)\\
			\sum_{v=M}^N Z_\zeta(t_{2v+1}) &= - 2N + O\big(N^{\frac{3}{4}} \log(N)^\frac{3}{4}\big).
		\end{align*}
	It follows that \( Z_\zeta(t_{2v}) \) is infinitely often positive while \( Z_\zeta(t_{2w+1}) \) is infinitely often negative. Hence there are infinitely many intervals \( (t_{2v}, t_{2w+1}] \), that contain an odd number of zeros of \( Z_\zeta(t) \) (counted with multiplicities). Since \( (t_{2v}, t_{2w+1}] \) is partitioned by an even number of intervals between consecutive Gram points, there are infinitely many intervals of the form \( (t_v, t_{v+1}] \) that contain an odd number of zeros of \( Z_\zeta(t) \). This fact is called the weak Gram law in some literature. It implies in particular that there are infinitely many zeros on the critical line, which has been proven first by Hardy \cite{hardy1}. For a survey on results regarding Gram's law we refer to Trudgian \cite{trudgian}.

	Our goal is to generalize this theorem and thereby the weak Gram law to the Hecke \( L \)-functions. These originate from modular forms and have properties similar to the Riemann zeta-function. In particular it is conjectured that an analogue of the Riemann hypothesis holds for them. We want to mention here that Lekkerkerker \cite{lekkerkerker} already studied the zeros of Hecke \( L \)-functions and showed among other things that there are infinitely many zeros on the respective critical line. For a given Hecke \( L \)-function \( L(s) \) we can define the analogous functions \( \vartheta_L(t), Z_L(t) \) and the corresponding Gram points \( t_v \) for \( v \geq v_0 \) (depending on \( L(s) \); see Section \ref{preparation} for the exact definitions) without difficulty.
	
	Now the result by Titchmarsh is proved with the aid of the approximate functional equation of \( \zeta(s) \) due to Hardy \& Littlewood \cite{hardy2}. An improvement of the error term in the theorem was recently achieved by Cao, Tanigawa and Zhai \cite{cao} with the help of a modified approximate functional equation with smooth weights. While approximate functional equations for Hecke \( L \)-functions have been proven for example by Apostol \& Sklar \cite{apostol} and Jutila \cite{jutila}, the error terms are not sufficiently small to use them in generalizing Titchmarsh's result. We therefore use a different approach that makes use of contour integration and leads to the following main theorem.
	
	\begin{theorem}\label{maintheorem} 
		We have for any \( \eps > 0 \), as \( T \to \infty \),
		\begin{align*}
			\sum_{T < t_{2v} \leq 2T} \omega(t_{2v}) Z_L(t_{2v}) &= \frac{1}{\pi} T + O_{L, \eps}\big(T^{\frac{3}{4}+\eps}\big)\\
			\sum_{T < t_{2v+1} \leq 2T} \omega(t_{2v+1}) Z_L(t_{2v+1}) &= -\frac{1}{\pi} T + O_{L, \eps}\big(T^{\frac{3}{4}+\eps}\big)\\
		\end{align*}
		with the weight function
		\[
			\omega(t) := \log \Big(\frac{t}{2\pi}\Big)^{-1}.
		\]
	\end{theorem}
	This theorem is sort of a weighted version of Titchmarsh's result for Hecke \( L \)-functions. By partial summation we can easily deduce an unweighted version from it.
	
	\begin{corollary}\label{maincorollary}
		We have for any fixed integer \( M \geq v_0 / 2 \) and any \( \eps > 0 \), as \( N \to \infty \),
		\begin{align*}
			\sum_{v=M}^N Z_L(t_{2v}) &= 2N + O_{L, \eps}\big(N^{\frac{3}{4}+\eps}\big)\\
			\sum_{v=M}^N Z_L(t_{2v+1}) &= - 2N + O_{L, \eps}\big(N^{\frac{3}{4}+\eps}\big).
		\end{align*}
	\end{corollary}
	From Theorem \ref{maintheorem} resp. Corollary \ref{maincorollary} the weak Gram law for Hecke \( L \)-functions follows as in the case of the Riemann zeta-function.
	
	Concerning notation in the following sections, \( \eps \) always denotes an arbitrarily small constant greater than \( 0 \), not necessarily the same at every occurrence. We write \( \int_z^w f(s) ds \) for the integral of \( f(s) \) along the the straight line form \( z \in \C \) to \( w \in \C \). Also we ommit the dependence of implicit and explicit constants on \( L \) and \( \eps \) for clarity.

\section{Preparation and preliminary results}\label{preparation}
	Let \( f(\tau) \) be a cusp form of weight \( k \geq 12 \) for the full modular group \( \SL_2(\Z) \) with the Fourier expansion
	\[
		f(\tau) = \sum_{n=1}^\infty a(n) e^{2\pi i n \tau},
	\]
	which additionally is a simultaneous eigenform of the Hecke-Operators. Then the coefficients \( a(n) \) are real and fulfil the bound \( a(n) = O(n^{\frac{k-1}{2} + \eps}) \) by Deligne \cite{deligne}. The corresponding Hecke \( L \)-function
	\[
		L(s) := \sum_{n=1}^\infty a(n) n^{-s}
	\]
	is absolutely convergent on the right half-plane \( \sigma > \frac{k+1}{2} \). It has an analytic continuation to the whole complex plane without poles and fulfils the functional equation
	\begin{equation}\label{FE}
		L(s) = \Delta_L(s) L(k-s) \quad \;\; \text{with} \;\; \Delta_L(s) := i^k (2\pi)^{2s -k} \frac{\Gamma(k-s)}{\Gamma(s)}.
	\end{equation}
	Hence the vertical line with real part \( \frac{k}{2} \) is the critical line of \( L(s) \). Also \( L(s) \) is a function of finite order on every vertical strip \( \sigma \in [\sigma_1, \sigma_2] \). For the theory of Hecke \( L \)-functions we refer to the monograph by Iwaniec \cite[Chapter 7]{iwaniec}.

	We now want to construct the analogues of the functions \( \vartheta_\zeta(t), Z_\zeta(t) \) and the Gram points \( t_v \) for the Hecke \( L \)-function \( L(s) \). For reasons that will become apparent later we do this with the help of a holomorphic logarithm of \( \Delta_L(s) \). We define this holomorphic logarithm explicitly with the unique holomorphic logarithm \( \log \Gamma(s) \), which is real-valued for real \( s \), by
	\begin{equation}\label{DefLogD}
		\log \Delta_L(s) := k \frac{\pi i}{2} + (2s - k) \log(2\pi) + \log \Gamma(k-s) - \log \Gamma(s)
	\end{equation}
	on the vertical strip \( \sigma \in (0, k) \). Now we can define the function \( \vartheta_L(t) \) for \( t \in \R \) by
	\begin{equation}\label{DefTheta}
		\vartheta_L(t) := \frac{i}{2} \log \Delta_L\Big(\frac{k}{2}+it\Big).
	\end{equation}
	From (\ref{FE}) we have \(|\Delta_L(\frac{k}{2} + it)| = 1\), hence \( \vartheta_L(t) \) is real-valued. Also by writing \( \Delta_L(s)^z := \exp(z \log\Delta_L(s)) \) for \( z \in \C \) we have
	\[
		e^{i \vartheta(t)} = \Delta_L\Big(\frac{k}{2} + it\Big)^{-\frac{1}{2}},
	\]
	so \( \vartheta_L(t) \) is a continuous branch of the argument of the function  \( \Delta_L(\frac{k}{2} + it)^{-\frac{1}{2}} \) analogously to \( \vartheta_\zeta(t) \). Now we define the continuous function
	\[
		Z_L(t) := e^{i \vartheta_L(t)} L\Big(\frac{k}{2}+it\Big)
	\]
	for \( t \in \R \). From \( \overline{L(s)} = L(\overline{s}) \) and the functional equation (\ref{FE}) it follows that \( Z_L(t) \) is also real-valued.
	
	To define the Gram points rigorously we need to show first that \( \vartheta_L(t) \) increases monotonically to infinity for \( t \) large enough. From Stirling's formula we can deduce the following approximations for \( \log \Delta_L(s) \) and its derivative.
	
	\begin{lemma}\label{LogDApprox}
		We have uniformly in the vertical strip \( \sigma \in (0, k) \), as \( t \to \infty \),
		\begin{align*}
			\log \Delta_L(s) &= (k - 2 \sigma) \log\Big(\frac{t}{2\pi}\Big) - 2 i t \log \Big( \frac{t}{2\pi e} \Big) + \frac{\pi i }{2} + O\Big(\frac{1}{t}\Big),\\
			\frac{d}{ds} \log \Delta_L(s) &= -2\log\Big(\frac{t}{2\pi}\Big) - \frac{i(k - 2\sigma)}{t} +  O\Big(\frac{1}{t^2}\Big).
		\end{align*}
		\begin{proof}
		From Stirlings formula (see \cite[Chapter XV, \S2, pp. 422-430]{lang}) we have uniformly in \( \sigma \in (0, k) \), as \( t \to \infty \),
		\[
			 \log \Gamma(s) = \Big(\sigma - \frac{1}{2}\Big) \log(t) - t\frac{\pi}{2} + \frac{1}{2} \log(2\pi) + i t \log\Big(\frac{t}{e}\Big) + i \Big( \sigma - \frac{1}{2}\Big) \frac{\pi}{2} + O\Big(\frac{1}{t}\Big)
		\]
		and
		\[
			 \frac{d}{ds} \log \Gamma(s) = \log(t) + \frac{\pi i}{2} - \frac{i( \sigma - \frac{1}{2} )}{t}  + O\Big(\frac{1}{t^2}\Big).
		\]
		Using this in (\ref{DefLogD}) and its derivate yields the approximations of \( \log \Delta_L(s) \) and \( \frac{d}{ds} \log \Delta_L(s) \) after lengthy computations.
		\end{proof}
	\end{lemma}
	
	By (\ref{DefTheta}) and its derivative
	\[
		\vartheta'_L(t) =\frac{d}{dt} \Big( \frac{i}{2} \log \Delta_L\Big(\frac{k}{2} + it\Big) \Big) = -\frac{1}{2} \frac{d}{ds} \log \Delta_L\Big(\frac{k}{2} + it\Big)
	\]
	we obtain the following corollary.
	
	\begin{corollary}\label{CorTheta}
		We have as \( t \to \infty \)
		\begin{align*}
			\vartheta_L(t) &= t \log \Big( \frac{t}{2\pi e} \Big) - \frac{\pi}{4} + O\Big(\frac{1}{t}\Big),\\
			\vartheta_L'(t) &= \log \Big( \frac{t}{2\pi} \Big) + O\Big(\frac{1}{t^2}\Big).
		\end{align*}
	\end{corollary}

	Hence the function \( \vartheta_L(t) \) increases monotonically for sufficiently large \( t \) and takes arbitrarily large values. We can therefore define the Gram points \( t_v \) of \( L(s) \) as the unique solutions of \( \vartheta_L(t_v) = v \pi \) for integers \( v \geq v_0 \) with some constant \( v_0 \in \Z \). We forgo a distinction to the Gram points of \( \zeta(s) \) in the notation for the sake of readability. We also need approximations for the Gram points \( t_v \) of \( L(s) \), their difference \( t_{v+1} - t_v \) and their number \( N(T) \) less than \( T \). These follow easily from corollary \ref{CorTheta} (analogous approximations for the Gram points of \( \zeta(s) \) are proven in \cite[\S6.1]{ivic1}).
	
	\begin{lemma}\label{LemmaGram}
		We have as \( v \to \infty \) resp. \( T \to \infty \)
		\[
			t_v \sim \frac{v \pi}{\log(v)}, \quad\;\; t_{v+1} - t_v \sim \frac{\pi}{\log(v)}, \quad\;\; N(T) \sim \frac{T\log(T)}{\pi}.
		\]
	\end{lemma}
	
	Now the idea of the proof of Theorem \ref{maintheorem} is to construct an auxiliary function \( G_L(s) \) with poles at \( \frac{k}{2} + i t_{2v} \), so that we can represent the sum of \( Z_L(t) \) at the Gram points \( t_{2v} \) as a contour integral by
	\[
		\res_{s = \frac{k}{2} + it_{2v}} \Big( G_L(s) \Delta_L(s)^{-\frac{1}{2}} L(s) \Big) = Z_L(t_{2v}) \res_{s = \frac{k}{2} + it_{2v}} G_L(s).
	\]
	Using the holomorphic logarithm \( \log \Delta_L(s) \) we define the auxiliary function as
	\[
		G_L(s) := - \frac{i}{2} \cot \Big(\frac{i}{4} \log \Delta_L(s) \Big),
	\]
	which is meromorphic on the vertical strip \( \sigma \in (0, k) \). On the critical line we have by (\ref{DefTheta})
	\begin{equation}\label{GCritical}
		G_L\Big(\frac{k}{2} + it\Big) =  - \frac{i}{2} \cot \Big(\frac{1}{2} \vartheta_L(t) \Big).
	\end{equation}
	\begin{lemma}\label{GPoles}
		For some constant \( A > 0 \) the poles of \( G_L(s) \) in the half-strip \( \sigma \in (0, k) \) and \( t > A \) lie exactly at \( s = \frac{k}{2} + i t_{2v} \) with \( t_{2v} > A \). For the residues we have
		\[
			\res_{s = \frac{k}{2} + it_{2v}} G_L(s) = \frac{1}{\vartheta'_L(t_{2v})} = \omega(t_{2v}) + O\Big( \frac{1}{t_{2v}^2} \Big).
		\]
	\begin{proof}
		Since the poles of the cotangent lie on the real axis, all poles \(s \) of \( G_L(s) \) fulfil \(\Real \log \Delta_L(s) = 0 \). From (\ref{FE}) we have \( |\Delta_L(\frac{k}{2} + it)| = 1 \), hence \( \Real \log \Delta_L(s) = 0 \) on the critical line. Furthermore Lemma \ref{LogDApprox} implies uniformly in \( \sigma \in (0, k) \)
		\[
			\frac{d}{d\sigma} \Real \log \Delta_L(\sigma + it) = \Real \frac{d}{ds} \log \Delta_L(\sigma + it) = - 2 \log\Big(\frac{t}{2\pi}\Big) + O\Big(\frac{1}{t^2}\Big).
		\] 
		Hence the function \( \frac{d}{d\sigma} \Real \log \Delta_L(\sigma + it) \) decreases monotonically with respect to \( \sigma \in (0, k) \) for fixed \( t > A \) with \( A \) being sufficiently large. Thus all the poles of \( G_L(s) \) for \( t > A \) lie on the critical line. From (\ref{GCritical}) it follows that the ordinates of these poles are exactly the Gram points with even index \( t_{2v} > A \).
		
		The poles \( \frac{k}{2} + it_{2v}\) for \( t > A \) are simple, since \( \vartheta_L(t) \) is increasing monotonically by Corollary \ref{CorTheta}, again assuming \( A \) to be sufficiently large. Hence we calculate for \( s = \frac{k}{2} + it \)
		\[
			\frac{d}{ds} \sin \Big( \frac{i}{4} \log \Delta_L(s) \Big) = (-i) \frac{d}{dt} \sin\Big(\frac{1}{2}\vartheta_L(t)\Big) = -\frac{i}{2} \cos\Big(\frac{1}{2}\vartheta_L(t)\Big) \vartheta'_L(t)
		\]
		and conclude
		\[
			\res_{s = \frac{k}{2} + i t_{2v}} G_L(s) = - \frac{i}{2} \cdot \frac{\cos\big(\frac{1}{2}\vartheta_L(t_{2v})\big)}{ -\frac{i}{2} \cos\big(\frac{1}{2}\vartheta_L(t_{2v})\big) \vartheta'_L(t_{2v})} = \frac{1}{\vartheta'_L(t_{2v})}.
		\]
		Again by Corollary \ref{CorTheta} we have
		\[
			\frac{1}{\vartheta'_L(t)} = \log\Big(\frac{t}{2\pi}\Big)^{-1} + O\Big( \frac{1}{t^2} \Big) = \omega(t) + O\Big( \frac{1}{t^2} \Big),
		\]
		from which the approximation of the residues follows.
	\end{proof}
	\end{lemma}
	
	At last we need a estimate for \( \Delta_L(s)^{-\frac{1}{2}} L(s) \), which follows from the Phragm\'en-Lindelöf principle (see \cite[Chapter XII, \S6]{lang}).
	
	\begin{lemma}\label{LemmaConvex}
		Let \( \frac{1}{2} < c < \frac{k}{2} \). Then we have 
		\[
			\Delta_L(s)^{-\frac{1}{2}} L(s) = O(t^c)
		\]
		uniformly in the strip \( \sigma \in \big[\frac{k}{2} - c, \frac{k}{2} + c\big] \) as \( t \to \infty \). In particular \( Z_L(t) = O\big(t^{\frac{1}{2} + \eps}\big)\).
	\begin{proof}
		Lemma \ref{LogDApprox} implies that, as \( t \to \infty \),
		\[
			\Delta_L\Big(\frac{k}{2} + c + it\Big)^{-\frac{1}{2}} = O\big(|t|^c\big).
		\]
		However, since \( \overline{\log \Delta_L(s)} = - k \pi i + \log \Delta_L(\overline{s}) \), which follows from (\ref{DefLogD}), this actually holds as \( |t| \to \infty \). Since the Dirichlet series of \( L(s) \) is absolutely convergent on the vertical line \( \sigma = \frac{k}{2} + c \) and hence bounded, we have
		\[
			\Delta_L\Big(\frac{k}{2} + c + it\Big)^{-\frac{1}{2}} L\Big(\frac{k}{2} + c + it\Big) = O\big(|t|^c\big)
		\]
		as \( |t| \to \infty \). The function \( \Delta_L(s)^{-\frac{1}{2}} L(s) \) takes the values of \( Z_L(t) \) on the critical line and thus is real-valued there. By the Schwarz reflection principle this yields additionally
		\[
			\Delta_L\Big(\frac{k}{2} - c + it\Big)^{-\frac{1}{2}} L\Big(\frac{k}{2} - c + it\Big) = O\big(|t|^c\big)
		\]
		as \( |t| \to \infty \). Since \( L(s) \) and \( \Delta_L(s)^{-\frac{1}{2}} \) are functions of finite order, we can apply the Phragm\'en-Lindelöf principle to \( \Delta_L(s)^{-\frac{1}{2}} L(s) \) in the vertical strip \( \sigma \in \big[\frac{k}{2} - c, \frac{k}{2} + c\big] \) and obtain \( \Delta_L(s)^{-\frac{1}{2}} L(s) = O(t^c) \) uniformly in this strip as \( t \to \infty \).
	\end{proof}
	\end{lemma}

\section{Proof of Theorem \ref{maintheorem}}

	Since we want to show an approximation for \( T \to \infty \) we can always assume that \( T > 0 \) is sufficiently large. Let \( T_0 \) and \( T_1 \) be Gram points with odd index, such that the intervals \( (T_0, T_1) \) and \( (T, 2T] \) contain the same Gram points with even index. In view of Lemma \ref{LemmaGram} we have
	\begin{equation}\label{T0T1}
		T_0 = T + O(1), \quad T_1 = 2T + O(1).
	\end{equation}
	Let \( \frac{1}{2} < c < \frac{k}{2} \) be a constant. We want to integrate the function \( G_L(s) \Delta_L(s)^{-\frac{1}{2}} L(s) \) along the positively oriented boundary of the rectangle \( \mathcal{R} \) with the vertices \( \frac{k}{2} \pm c + i T_0 \) and \( \frac{k}{2} \pm c + i T_1 \). By Cauchy's residue theorem we obtain in view of Lemma \ref{GPoles}
		\begin{equation}\label{Integral1}
			\sum_{T < t_{2v} \leq 2T} \res_{s = \frac{k}{2} + i t_{2v}} \Big(  G_L(s) \Delta_L(s)^{-\frac{1}{2}} L(s)  \Big) = \frac{1}{2\pi i} \int_{\partial \mathcal{R}} G_L(s) \Delta_L(s)^{-\frac{1}{2}} L(s) ds.
		\end{equation}
		We first deal with the left-hand side. Lemma \ref{GPoles} gives
		\[
			\res_{s = \frac{k}{2} + i t_{2v}} \Big(  G_L(s) \Delta_L(s)^{-\frac{1}{2}} L(s)  \Big) = \omega(t_{2v}) Z_L(t_{2v}) + O\Bigg(\frac{Z_L(t_{2v})}{t_{2v}^2}\Bigg). 
		\]
		Using Lemma \ref{LemmaConvex} and Lemma \ref{LemmaGram} we obtain for the sum of the error terms 
		\[
			\sum_{T < t_{2v} \leq 2T} \frac{|Z_L(t_{2v})|}{t_{2v}^2} \ll \sum_{T < t_{2v} \leq 2T} t_{2v}^{-\frac{3}{2} + \eps} \ll T^{-\frac{3}{2} + \eps} N(2T) \ll T^{-\frac{1}{2} + \eps}.
		\]
		Hence the left-hand side of (\ref{Integral1}) is
		\[
			\sum_{T < t_{2v} \leq 2T} \omega(t_{2v}) Z_L(t_{2v}) + O\big(T^{-\frac{1}{2} + \eps}\big).
		\]
		On the right-hand side of (\ref{Integral1}) we split the integral into the four integrals along the sides of \( \mathcal{R} \). First we want to estimate the integrals along the horizontal sides, which have the form
		\[
			\int_{\frac{k}{2} - c + it}^{\frac{k}{2} + c + it} G_L(s) \Delta_L(s)^{-\frac{1}{2}} L(s) ds
		\]
		with \( t = t_{2v+1} \) for an integer \( v \geq v_0 \). For that we show that the function \( G_L(s) \) is bounded on the horizontal paths \( \sigma \mapsto \sigma + i t_{2v+1} \) with \( \sigma \in [\frac{k}{2} - c, \frac{k}{2} + c] \) as \( t_{2v+1} \to \infty \). For \( \sigma = \frac{k}{2} \) we have
		\[
			\frac{i}{4} \log \Delta_L\Big(\frac{k}{2} + it_{2v+1}\Big) = \frac{1}{2} \theta_L(t_{2v+1}) = \Big(v+\frac{1}{2}\Big) \pi
		\]
		and the real part of \( \frac{i}{4} \log \Delta_L(s) \) is independent of \( \sigma \) except for the error term \( O(t^{-1}) \) by Lemma \ref{LogDApprox}. Hence the real part of \( \frac{i}{4} \log \Delta_L(\sigma + it_{2v+1}) \) for \( \sigma \in [\frac{k}{2} - c, \frac{k}{2} + c] \) lies in the interval \( ( v +  [\frac{1}{4}, \frac{3}{4}]) \pi \) for sufficiently large \( t_{2v+1} \). The cotangent \( \cot(x + i y) \) is bounded in the vertical strips \( x \in ( v +  [\frac{1}{4}, \frac{3}{4}]) \pi  \) because of its periodicity and 
		\[
			|\cot(x+iy)|= \frac{|e^{i(x+iy)} + e^{-i(x+iy)}|}{|e^{i(x+iy)} - e^{-i(x+iy)}|} \leq \frac{e^{y} + e^{-y}}{|e^y - e^{-y}|} \ll \frac{e^{|y|}}{ e^{|y|}} = 1
		\]
		as \( y \to \infty \). Therefore we have \( G_L(s) = O(1) \) on the horizontal paths as \( t_{2v+1} \to \infty \). By Lemma \ref{LemmaConvex} it follows
		\[
			\int_{\frac{k}{2} - c + it}^{\frac{k}{2} + c + it} G_L(s) \Delta_L(s)^{-\frac{1}{2}} L(s) ds = O(t^c)
		\]
		as \( t = t_{2v+1} \to \infty \). In view of (\ref{T0T1}) the horizontal integrals on the right-hand side of (\ref{Integral1}) are therefore bounded by \( O(T^c) \).

		Next we deal with the integral along the left vertical side of \( \mathcal{R} \). Using the functional equation (\ref{FE}) and \( \overline{L(s)} = L(\overline{s}) \) we obtain
		\begin{equation}\label{LeftIntegral}
			\int_{\frac{k}{2} - c + iT_1}^{\frac{k}{2} - c + iT_0} G_L(s) \Delta_L(s)^{-\frac{1}{2}} L(s) ds = \int_{\frac{k}{2} - c + iT_1}^{\frac{k}{2} - c + iT_0} G_L(s) \Delta_L(s)^{\frac{1}{2}} \overline{L(k-\overline{s})} ds.
		\end{equation}
		From (\ref{DefLogD}) we have \( \log \Delta_L(s) = - \overline{\log \Delta_L(k-\overline{s})} \). Using this we easily obtain the functional equations \( G_L(s) = -\overline{G_L(k - \overline{s})} \) and \(\Delta_L(s)^{1/2} = \overline{ \Delta_L(k - \overline{s})^{-1/2}} \). Together with a parametrization \( s = \frac{k}{2} - c + it \) of the path of integration we obtain through delicate transformations
		\[
			\int_{\frac{k}{2} - c + iT_1}^{\frac{k}{2} - c + iT_0} G_L(s) \Delta_L(s)^{\frac{1}{2}} \overline{L(k-\overline{s})} ds = - \overline{\int_{\frac{k}{2} + c + iT_0}^{\frac{k}{2} + c + iT_1} G_L(s) \Delta_L(s)^{-\frac{1}{2}} L(s) ds}.
		\]	
		Hence the left vertical integral is equal to the negative conjugate of the right vertical integral. Altogether we have transformed (\ref{Integral1}) to
		\begin{equation}\label{Integral2}
			\sum_{T < t_{2v} \leq 2T} \omega(t_{2v}) Z_L(t_{2v}) = \frac{1}{\pi} \Imag \Big( \int_{\frac{k}{2} + c + iT_0}^{\frac{k}{2} + c + iT_1} G_L(s) \Delta_L(s)^{-\frac{1}{2}} L(s) ds \Big) + O(T^c).
		\end{equation}
		Now we approximate the term \( G_L(s) \Delta_L(s)^{-\frac{1}{2}} \) in the integrand for \( s = \frac{k}{2} + c + it \) and \( t \to \infty \). Substituting \( z := \frac{i}{4} \log\Delta_L(s) \) gives
		\[
			G_L(s) \Delta_L(s)^{-\frac{1}{2}} 
			= - \frac{i}{2} \cot(z) e^{2iz}
			= \frac{1}{2} \cdot \frac{e^{2iz} +1}{e^{2iz}-1} e^{2iz}
			= \frac{1}{2} e^{2iz} + 1 + \frac{1}{e^{2iz}-1}.
		\]
		By Lemma \ref{LogDApprox} the imaginary part \( y = \Imag(z) \) for \( s = \frac{k}{2} + c + it \) is
		\[
			y = \frac{1}{4} \Real \log\Delta_L(s) = -\frac{c}{2} \log\Big(\frac{t}{2\pi}\Big) + O\Big(\frac{1}{t}\Big).
		\]
		Hence \( y \) grows in the negative direction and \( |e^{2iz} - 1| \geq |e^{-2y} - 1| \gg e^{-2y} \) as \( t \to \infty \). Thus we obtain the approximation
		\begin{equation}\label{GApprox}
			G_L(s) \Delta_L(s)^{-\frac{1}{2}} = \frac{1}{2} e^{2iz} + 1 + O(e^{2y}) = \frac{1}{2} \Delta_L(s)^{-\frac{1}{2}} + 1 + O(t^{-c})
		\end{equation}
		as \( t \to \infty \). Using this in (\ref{Integral2}) yields
		\begin{equation}\label{Integral3}
		\begin{split}
			\sum_{T < t_{2v} \leq 2T} \omega(t_{2v}) Z_L(t_{2v}) 
			&= \frac{1}{2 \pi} \Imag \Big( \int_{\frac{k}{2} + c + iT_0}^{\frac{k}{2} + c + iT_1} \Delta_L(s)^{-\frac{1}{2}} L(s) ds \Big) \\
			&+ \frac{1}{\pi} \Imag \Big( \int_{\frac{k}{2} + c + iT_0}^{\frac{k}{2} + c + iT_1} L(s) ds \Big) \\
			&+ O(T^c).
		\end{split}
		\end{equation}
		Here we have used that \( L(s) \) is bounded on the vertical line \( \sigma = \frac{k}{2} + c \) because of absolute convergence. We compute the second integral in (\ref{Integral3}) using the Dirichlet series \( L(s) = \sum_{n=1}^\infty a(n) n^{-s} \) with \( a(1) = 1 \). Interchanging integration and summation by the theorem of Lebesgue then yields
		\begin{align*}
			\frac{1}{\pi} \Imag \Big( \int_{\frac{k}{2} + c + iT_0}^{\frac{k}{2} + c + iT_1} L(s) ds \Big)
			&= \frac{1}{\pi} \Real \Big( \int_{T_0}^{T_1} \sum_{n=1}^\infty a(n) n^{-\frac{k}{2} - c - it} dt \Big)\\
			&= \frac{1}{\pi} \Real \Big( \sum_{n=1}^\infty a(n) n^{-\frac{k}{2} - c} \int_{T_0}^{T_1} n^{-it} dt \Big) \\
			&= \frac{1}{\pi} (T_1 - T_0) + \Real \Big( \sum_{n=2}^\infty a(n) n^{-\frac{k}{2} - c} O(1) \Big).
		\end{align*}
		Using \( T_1 - T_0 = T + O(1) \) by (\ref{T0T1}) and the absolute convergence of \( L(s) \) at \( s = \frac{k}{2} + c \) we obtain
		\[
			\frac{1}{\pi} \Imag \Big( \int_{\frac{k}{2} + c + iT_0}^{\frac{k}{2} + c + iT_1} L(s) ds \Big)
			= \frac{1}{\pi} T + O(1).
		\]
		We also want to interchange \( T_0 \) with \( T \) and \( T_1 \) with \( 2 T \) in the first integral of (\ref{Integral3}). Since both differences lie in \( O(1) \) by (\ref{T0T1}) and \( \Delta_L(s)^{-\frac{1}{2}} L(s) = O(t^c) \) by Lemma \ref{LemmaConvex}, this yields again the error term \( O(T^c) \). Hence we have transformed (\ref{Integral3}) to
		\begin{equation}\label{Integral4}
			\sum_{T < t_{2v} \leq 2T} \omega(t_{2v}) Z_L(t_{2v}) 
			= \frac{1}{\pi} T + \frac{1}{2 \pi} \Imag \Big( \int_{\frac{k}{2} + c + iT}^{\frac{k}{2} + c + 2iT} \Delta_L(s)^{-\frac{1}{2}} L(s) ds \Big)  + O(T^c).
		\end{equation}
		It remains to estimate the integral of \( \Delta_L(s)^{-\frac{1}{2}} L(s) \). By Lemma \ref{LogDApprox} we have
	\[
			\Delta_L\Big(\frac{k}{2} + c + it \Big)^{-\frac{1}{2}}
			= e^{-\frac{\pi i}{4}} \Big( \frac{t}{2\pi} \Big)^c \exp \Big( i t \log \Big( \frac{t}{2\pi e} \Big) \Big) \Big( 1 + O \Big(\frac{1}{t}\Big) \Big).
	\]
	We use this approximation in the integral of \( \Delta_L(s)^{-\frac{1}{2}} L(s) \) and proceed as in the estimation of the second integral of (\ref{Integral3}). This yields
	\begin{equation}\label{Series}
		 \int_{\frac{k}{2} + c + i T}^{\frac{k}{2} + c + 2 i T} \Delta_L(s)^{-\frac{1}{2}} L(s) ds = 2 \pi e^\frac{\pi i}{4} \sum_{n=1}^\infty a(n) n^{-\frac{k}{2}-c} I(n) + O(T^c)
	\end{equation}
	with
	\begin{equation}\label{DefI}
		I(n) := \frac{1}{2\pi} \int_{T}^{2T} \Big(\frac{t}{2\pi}\Big)^c \exp \Big( i t \log \Big( \frac{t}{2 \pi e n} \Big) \Big) dt.
	\end{equation}
	We also define \( \hat{T} := \frac{T}{2\pi} \) and \( F_n(t) := t \log \big( \frac{t}{e n} \big) \). By a change of variables we can then rewrite (\ref{DefI}) as
	\[
		I(n) = \int_{\hat{T}}^{2\hat{T}} t^c \exp( 2 \pi i F_n(t) ) dt.
	\]
	The function \( F_n'(t) = \log(\frac{t}{n}) \) has the zero \( t = n \), which lies in the interval of integration for \( \hat{T} \leq n \leq 2 \hat{T} \). In view of this saddle point we split the series on the right-hand side of (\ref{Series}) into
	\begin{equation}\label{SeriesSplit}
		\sum_{n=1}^\infty a(n) n^{-\frac{k}{2}-c} I(n) = \sum\nolimits_1 + \sum\nolimits_2 + \sum\nolimits_3 + \sum\nolimits_4 + \sum\nolimits_5,
	\end{equation}
	where the ranges of summation, depending on a constant \( d \in (0, 1) \), are the following:
	\begin{alignat*}{5}
		\sum\nolimits_1 &: &\quad 1 &\leq n & &\leq  \hat{T} - \hat{T}^d, \\
		\sum\nolimits_2 &: &\quad  \hat{T} - \hat{T}^d &< n & &\leq  \hat{T} + \hat{T}^d, \\
		\sum\nolimits_3 &: &\quad  \hat{T} + \hat{T}^d &< n & &\leq  2 \hat{T} - \hat{T}^d, \\
		\sum\nolimits_4 &: &\quad  2 \hat{T} - \hat{T}^d &< n & &\leq  2 \hat{T} + \hat{T}^d, \\
		\sum\nolimits_5 &: &\quad  2 \hat{T} + \hat{T}^d &< n. & & 
	\end{alignat*}
	First let \( 1 \leq n \leq \hat{T} - \hat{T}^d \). Then the function \( F_n'(t) = \log(\frac{t}{n}) \) grows monotonically in the range \( t \in [\hat{T}, 2\hat{T}] \) and fulfils
	\[
	 	F_n'(t) = \log\Big(\frac{t}{n}\Big) \geq \log\Big(\frac{\hat{T}}{\hat{T}-\hat{T}^d}\Big) = - \log(1 - \hat{T}^{d-1}) \asymp \hat{T}^{d-1}.
	 \]
	 Applying the first derivative test (see \cite[Lemma 2.1]{ivic2}) yields
	 \[
	 	 I(n) = \int_{\hat{T}}^{2\hat{T}} t^c \exp( 2 \pi i F_n(t) ) dt \ll \hat{T}^{c + 1 - d}.
	 \]
	 Thus we obtain using \( a(n) = O(n^{\frac{k-1}{2}+\eps}) \)
	 \[
	 	\sum\nolimits_1 = \sum_{1 \leq n \leq \hat{T} - \hat{T}^d} a(n) n^{-\frac{k}{2}-c} I(n) \ll \hat{T}^{c + 1 - d} \sum_{n = 1}^\infty n^{-\frac{1}{2} - c + \eps} \ll \hat{T}^{c + 1 - d}.
	 \]
	 Hence \( \sum_1 \in O(T^{c+1-d}) \) and in a similar way \( \sum_5 \in O(T^{c + 1 - d}) \) follows. 
	 
	 Now let \( \hat{T} - \hat{T}^d < n \leq  \hat{T} + \hat{T}^d \). Then \( F_n''(t) = t^{-1} \gg \hat{T}^{-1} \) and an application of the second derivate test (see \cite[Lemma 2.2]{ivic2}) yields
	\[
		I(n) = \int_{\hat{T}}^{2\hat{T}} t^c \exp( 2 \pi i F_n(t) ) dt  \ll \hat{T}^{c + \frac{1}{2}}.
	\]
	Again using \( a(n) = O(n^{\frac{k-1}{2}+\eps}) \) we obtain  
	\[
		\sum\nolimits_2 = \sum_{\hat{T} - \hat{T}^d < n \leq  \hat{T} + \hat{T}^d} a(n) n^{-\frac{k}{2}-c} I(n)
		\ll \hat{T}^{c + \frac{1}{2}} \sum_{\hat{T} - \hat{T}^d < n \leq  \hat{T} + \hat{T}^d} n^{-\frac{1}{2} - c + \eps}
		\ll \hat{T}^{d + \eps}.
	\]
	Hence \( \sum_2 = O(T^{d+\eps}) \) and \( \sum_4 = O(T^{d+\eps}) \) follows analogously.
	
	It remains to estimate \( \sum_3 \). We use the following lemma from \cite[Lemma III.\S1.2]{karatsuba}.
	
	\begin{lemma}\label{LemmaKaratsuba}
		Suppose that \( f(t) \) and \( \varphi(t) \) are real-valued functions on the interval \( [a, b] \) which satisfy the conditions
		\begin{enumerate}
			\item \( f^{(4)}(t) \) and \( \varphi''(t) \) are continuous,
			
			\item there exists \( 0 < b-a \leq U \), \( 0 < H \), \( A < U \), such that
				\begin{alignat*}{5}
					f''(t) &\asymp A^{-1}, &\quad\;\; f^{(3)}(t) &\ll A^{-1} U^{-1}, &\quad\;\; f^{(4)}(t) &\ll A^{-1} U^{-2}, \\
					g(t) & \ll H, &\quad\;\; g'(t) &\ll H U^{-1}, &\quad\;\; g''(t) &\ll H U^{-2},
				\end{alignat*}
			\item \( f'(t_0) = 0 \) for some \( t_0 \in [a, b] \).
		\end{enumerate}
		Then
		\begin{align*}
			\int_a^b \varphi(t) \exp(2 \pi i f(t)) dt 
			&= \frac{\varphi(t_0)}{\sqrt{f''(t_0)}} \exp \Big(2 \pi i f(t_0) + \frac{\pi i}{4}\Big) + O(H A U^{-1})\\
			&+ O\Big( H \cdot\min \big(|f'(a)|^{-1}, \sqrt{A} \, \big) \Big)
			+ O\Big( H \cdot \min \big(|f'(b)|^{-1}, \sqrt{A} \, \big) \Big).
		\end{align*}
	\end{lemma}

	Now let \(\hat{T} + \hat{T}^d < n \leq  2 \hat{T} - \hat{T}^d \). We have
	\begin{alignat*}{5}
		F_n(t) &= t \log\Big(\frac{t}{en}\Big), &\quad\;\; F_n'(t) &= \log\Big(\frac{t}{n}\Big), &\quad\;\;  &\quad \;\; \\
		F_n''(t) &= \frac{1}{t}, &\quad\;\; F_n^{(3)}(t) &= -\frac{1}{t^2}, &\quad\;\; F_n^{(4)}(t) &= \frac{2}{t^3}.		
	\end{alignat*} 
	Applying Lemma \ref{LemmaKaratsuba} with \( f(t) = F_n(t), \varphi(t) = t^c  \) and \( A = \hat{T}, U = 2\hat{T}, H=\hat{T}^c,  t_0 = n \) yields
	\begin{equation}\label{Integral5}
	\begin{split}
		I(n) &= n^{c + \frac{1}{2}} \exp\Big( - 2\pi i n + \frac{\pi i}{4} \Big) + O(\hat{T}^c)\\
		&+ O\Big( \hat{T}^c \cdot \min \big(|F_n'(\hat{T})|^{-1}, \sqrt{\hat{T}}\, \big) \Big)
		+ O\Big( \hat{T}^c \cdot \min \big(|F_n'(2\hat{T})|^{-1}, \sqrt{\hat{T}}\, \big) \Big).
	\end{split}
	\end{equation}
	We have
	\begin{align*}
		|F_n'(\hat{T})| &= \Big| \log\Big( \frac{\hat{T}}{n} \Big) \Big| = \log\Big( \frac{n}{\hat{T}} \Big) \geq \log(1 + \hat{T}^{d-1}) \asymp \hat{T}^{d-1},\\
		|F_n'(2\hat{T})| &= \log\Big( \frac{2\hat{T}}{n} \Big) \geq \log \Big( \frac{2\hat{T}}{2 \hat{T} - \hat{T}^d} \Big) = - \log\Big( 1 - \frac{1}{2}\hat{T}^{d-1} \Big) \asymp \hat{T}^{d-1},
	\end{align*}
	hence \( |F_n'(\hat{T})|^{-1}, |F_n'(2\hat{T})|^{-1} \ll \hat{T}^{1-d} \). This gives the overall error term \( O(\hat{T}^{c+1-d}) \) in (\ref{Integral5}), which is independent of \( n \). Therefore
	\begin{equation}\label{sum}
		\sum\nolimits_3 = e^{\frac{i\pi}{4}} \sum_{\hat{T} + \hat{T}^d < n \leq  2 \hat{T} - \hat{T}^d} a(n) n^{-\frac{k-1}{2}} + O( \hat{T}^{c+1-d}),
	\end{equation}
	where we have used the absolute convergence of \( L(s) \) at \( s = \frac{k}{2} + c \). 
	
	It remains to deal with the sum of \( a(n) n^{-\frac{k-1}{2}} \). We need to use a fact about the coefficients of cusp forms of weight \( k \), namely that
	\[
		\sum_{n \leq x} a(n) \ll x^\frac{k}{2} \log(x)
	\]
	as \( x \to \infty \) (see \cite[Theorem 5.3]{iwaniec}). By partial summation we then obtain for the sum on the right-hand side of (\ref{sum})
	\[
		\sum_{\hat{T} + \hat{T}^d < n \leq  2 \hat{T} - \hat{T}^d} a(n) n^{-\frac{k-1}{2}} \ll \hat{T}^{\frac{1}{2} + \eps},
	\]
	hence \( \sum_3 = O(T^{c+1-d}) \). From \( \sum_1, \sum_3, \sum_5 = O(T^{c + 1 - d}) \) and \( \sum_2, \sum_4 = O(T^{d+\eps}) \) it follows in view of (\ref{Series}) and (\ref{SeriesSplit}) that
	\[
		\int_{\frac{k}{2} + c + i T}^{\frac{k}{2} + c + 2 i T} \Delta_L(s)^{-\frac{1}{2}} L(s) ds = O(T^{c+1-d}) + O(T^{d+\eps}) + O(T^c).
	\]
	We choose \( c = \frac{1}{2} + \eps \) and \( d = \frac{3}{4} \) to obtain the overall bound \( O(T^{\frac{3}{4} + \eps}) \). Then (\ref{Integral4}) gives the final approximation
	\[
		\sum_{T < t_{2v} \leq 2T} \omega(t_{2v}) Z_L(t_{2v}) 
		= \frac{1}{\pi} T + O(T^{\frac{3}{4}+\eps}).
	\]
	Hence the treatment of the first sum in Theorem \ref{maintheorem} is finished. We can deal with the second sum analogously using the auxiliary function
	\[
		H(s) := \frac{i}{2} \tan \Big( \frac{i}{4} \log \Delta_L(s) \Big)
	\]
	instead of \( G(s) \). Then the approximation 
	\[
			H(s) \Delta_L(s)^{-\frac{1}{2}}= \frac{1}{2} \Delta_L(s)^{-\frac{1}{2}} - 1 + O(t^{-c})
	\]
	for \( s = \frac{k}{2} + c + it \) as \( t \to \infty \) in comparison with (\ref{GApprox}) leads to the negative dominant term in the approximation of the second sum.

\section{Proof of Corollary \ref{maincorollary}}
		We consider
		\[
			S(T) := \sum_{t_{2v} \leq T} \omega(t_{2v}) Z_L(t_{2v}),
		\]
		where the summation ranges over all Gram points \( t_{2v} \) less than or equal to \( T \). By Theorem \ref{maintheorem} we obtain the approximation
		\begin{equation}\label{ApproxS}
			S(T) 
			= \sum_{m=1}^\infty \; \sum_{T/2^m < t_{2v} \leq 2T/2^m} \omega(t_{2v}) Z_L(t_{2v}) 
			= \frac{T}{\pi} + O(T^{\frac{3}{4}+\eps}).
		\end{equation}
		Now we deal with the sum of \( Z_L(t_{2v}) \) for the Gram points \( t_{2v} \), where \( M \leq v \leq N \) resp. \( t_{2M} \leq t_{2v} \leq t_{2N} \). An application of partial summation yields
		\begin{align*}
			\sum_{v = M}^N Z_L(t_{2v}) 
			&= \sum_{t_{2M} \leq t_{2v} \leq t_{2N}} \log\Big(\frac{t_{2v}}{2\pi}\Big) \omega(t_{2v}) Z_L(t_{2v}) \\
			&= \log\Big(\frac{t_{2N}}{2\pi}\Big) S(t_{2N}) - \int_{0}^{t_{2N}} \frac{S(T)}{T} dT + O(1).
		\end{align*}
		Using (\ref{ApproxS}) and the estimate \( t_{2N} \ll N \), which follows from Lemma \ref{LemmaGram}, we obtain
		\begin{equation}\label{ApproxSum}
		\begin{split}
			\sum_{v = M}^N Z_L(t_{2v}) 
			&= \frac{t_{2N}}{\pi} \log\Big(\frac{t_{2N}}{2\pi}\Big) - \frac{ t_{2N}}{\pi} + O(N^{\frac{3}{4}+\eps}) \\
			&= \frac{t_{2N}}{\pi} \log\Big(\frac{t_{2N}}{2\pi e}\Big) + O(N^{\frac{3}{4}+\eps}).
		\end{split}
		\end{equation}
		By the definition of the Gram points and Corollary \ref{CorTheta} we have
		\[
			v = \frac{1}{\pi} \theta_L(t_v) = \frac{t_v}{\pi} \log\Big(\frac{t_v}{2\pi e}\Big) + O(1).
		\]
		Using this for the Gram point \( t_{2N} \) in (\ref{ApproxSum}) gives
		\[
			\sum_{v = M}^N Z_L(t_{2v}) = 2N + O(N^{\frac{3}{4}+\eps}).
		\]
		The approximation of the second sum follows analogously.

\section{Concluding remarks}	
		The most difficult part in the proof of Theorem \ref{maintheorem} is to estimate the integral in (\ref{Integral4}). By shifting the path of integration to the left onto the critical line using Lemma \ref{LemmaConvex} we have
		\begin{equation}\label{Integral4Z}
			\Imag \Big( \int_{\frac{k}{2} + c + i T}^{\frac{k}{2} + c + 2 i T} \Delta_L(s)^{-\frac{1}{2}} L(s) ds \Big) = \int_T^{2T} Z_L(t) dt + O(T^c).
		\end{equation}
		Ivi\'{c} \cite{ivic3} showed that the integral of \( Z_\zeta(t) \) over the interval \( [T, 2T] \) is bounded by \( O(T^{\frac{1}{4}+\eps})\). He mentioned a possible but insufficient approach in his article, which we have adopted to deal with the integral in (\ref{Integral4}). Hence generalizing the actual method by Ivi\'{c} might yield an improvement of the error term in Theorem \ref{maintheorem}. Also in view of (\ref{Integral4Z}) we have showed implicitly, that the integral of \( Z_L(t) \) over the interval \( [T, 2T] \) is bounded by \( O(T^{\frac{3}{4}+\eps})\).

		From Theorem 1 an analogous result for the Hecke \( L \)-functions \( L(s) \) of arbitrary cusp forms follows, since every cusp form of weight \( k \) for the full modular group \( \SL_2(\Z) \) is a linear combination of simultaneous eigenforms of the Hecke-Operators with complex coefficients. Then the dominant terms are \( \pm \frac{a(1)}{\pi}T \), where \( a(1) \) is the first coefficient of \( L(s) \). If the cusp form is a linear combination of simultaneous eigenforms with real coefficients, the analogously defined function \( Z_L(t) \) is also real-valued. Hence in this case the weak Gram law for \( L(s) \) follows, provided that \( a(1) \neq 0 \).

		Also Theorem \ref{maintheorem} can be generalized to Hecke \( L \)-functions corresponding to cusp forms for congruence subgroups of \( \SL_2(\Z) \) without difficulty.

	\printbibliography[heading=bibintoc]

\end{document}